\documentclass[12pt]{article}
\usepackage{amsmath, amssymb, amsthm, enumerate} 
\usepackage[utf8]{inputenc}
\usepackage{mathtools}

\newtheorem{theorem}{Theorem}
\newtheorem{proposition}[theorem]{Proposition}

\newtheorem{lemma}[theorem]{Lemma}

\def\rp{{\mathbb R}^p}

\def\cc{{\mathbb C}}
\def\zz{{\mathbb Z}}

\def\su{\subset}

\def\se{\setminus}

\def\de{\delta}
\def\De{\Delta}

\def\Si{\Sigma}

\def\om{\omega}

\def\cd{\cdot}
\def\stb{,\ldots ,}

\def\emp{\emptyset}

\def\ran{\rangle}
\def\lan{\langle}

\def\akkor{\Longrightarrow}

\def\ol{\overline}

\def\sumin{\sum_{i=1}^n}
\def\sumis{\sum_{i=1}^s}

\def\sumik{\sum_{i=1}^k}
\def\sumjk{\sum_{j=1}^k}
\def\sumi0n{\sum_{i=0}^n}

\def\x1n{x_1 \stb x_n}
\def\y1n{y_1 \stb y_n}
\def\deg{{\rm deg}\, }
\def\dim{{\rm dim}\, }

\def\det{{\rm det}\, }
\def\spect{{\rm sp}\, }

\begin{document}

\title{A characterization of generalized exponential polynomials in terms
  of decomposable functions}

\author{{\bf Mikl\'os Laczkovich} (Budapest, Hungary)}

\maketitle

\begin{abstract}
Let $G$ be a topological commutative semigroup with unit.
We prove that a continuous function $f\colon G\to \cc$ is a
generalized exponential polynomial if and only if there is an $n\ge 2$
such that $f(x_1 +\ldots +x_n )$ is decomposable; that is, if
$f(x_1 +\ldots +x_n )=\sumik u_i \cd v_i$, where the function
$u_i$ only depends on the
variables belonging to a set $\emp \ne E_i \subsetneq \{ x_1 \stb x_n \}$,
and $v_i$ only depends on the variables belonging to $\{ x_1 \stb x_n \} \se
E_i$ $(i=1\stb k)$.
\end{abstract}

\section{Introduction and main results}
Let $X$ be a nonempty set. A function of $n$ variables $F\colon X^n \to \cc$
is called {\it decomposable of order $k$} if there are functions
$u_i ,v_i \colon X^n \to \cc$ $(i=1\stb k)$ such that $F=\sumik u_i \cd v_i$ and,
for every $i=1\stb k$, $u_i$ only depends on the
variables belonging to a set $\emp \ne E_i \subsetneq \{ x_1 \stb x_n \}$,
and $v_i$ only depends on the variables belonging to $\{ x_1 \stb x_n \} \se
E_i$. For example, if $G$ is a semigroup and $f\colon G\to \cc$ 
satisfies a Levi-Civita equation $f(xy)=\sumik a_i (x)\cd b_i (y)$
$(x,y\in G)$, then $f(xy)$ is decomposable of order $k$.

The following remarkable result was proved by
Ekaterina Shulman in \cite[Theorem 6]{Sh}.
{\it Let $G$ be a topological
semigroup with unit. If a continuous function $f\colon G\to \cc$ is such that $f
(x_1 \cdots x_n )$ is decomposable for some $n>1$, then $f$ is an
almost matrix function.}

A function is a matrix function if $f(xy)$ satisfies a Levi-Civita
equation. The function $f$ is an almost matrix function
if, for every finite set $E\su G$, there is a finite dimensional subspace
of $C(G)$ containing $f$ and invariant under the subsemigroup generated by $E$.
It was also shown by Shulman that {\it if $G$ is the union of an increasing
net of
topologically $p$-generated subsemigroups, then, whenever $f\in C(G)$ is
such that $f(x_1 \cdots x_n )$ is decomposable for some $n>1$, then $f$ is a
matrix function.} In particular, {\it if $f\in C(\rp )$ is such that $f(x_1
+\ldots +x_n )$ is decomposable for some $n>1$, then $f$ is an
exponential polynomial} (\cite[Corollary 8]{Sh}). 

In this note our aim is to prove that in the commutative case, if
$f(x_1 +\ldots +x_n )$ is decomposable, then $f$ is necessarily a
generalized exponential polynomial. The relevant definitions are as follows.

Let $G$ be an Abelian topological semigroup with unit, written additively.
Let $C(G)$ denote the set of complex valued continuous functions defined on $G$.
The difference operator $\De _h$ is defined by $\De _h f (x)=f(x+h)-f(x)$
$(x\in G)$ for every $h\in G$ and $f\colon G\to \cc$. A function $f\in C(G)$
is a {\it generalized polynomial}, if there is an $n\ge 0$ such that
$\De _{h_1} \ldots \De _{h_{n+1}} f=0$ for every $h_1 \stb h_{n+1} \in G$. The
smallest
$n$ with this property is the {\it degree of $f$}, denoted by $\deg f$.
The degree of the identically zero function is $-1$.
The set of generalized polynomials will be denoted by \hbox{GP}.

A function $m\in C(G)$ is an {\it exponential}, if $m\ne 0$ and
$m(x+y)=m(x)\cd m(y)$ for every $x,y\in G$. The functions of the form
$\sumis p_i \cd m_i$, where $p_i \in \hbox{GP}$ and $m_i$ is an exponential
for every $i=1\stb s$ are called {\it generalized exponential polynomials.}
The set of generalized exponential polynomials will be denoted by \hbox{GEP}.

{\it In the rest of the paper we assume that $G$ is a commutative unital
semigroup written additively.} Our main result is the following.

\begin{theorem}\label{t1}
For every continuous function $f\in C(G)$ the following are equivalent.
\begin{enumerate}[{\rm (i)}]
\item There is an $n\ge 2$ such that $f(x_1 +\ldots +x_n )$ is decomposable.
\item $f$ is a generalized exponential polynomial.
\end{enumerate}
\end{theorem}
It is easy to see that if $f(x_1 +\ldots +x_n )$ is decomposable, then so is
$f(x_1 +\ldots +x_{n'} )$ for every $n'>n$. Thus (i) is equivalent to the
statement that $f(x_1 +\ldots +x_n )$ is decomposable for every $n$ large
enough.

Note that the class \hbox{GEP} is contained in the class of
almost matrix functions, and the containment, in general, is strict
(see \cite[example 1, p. 18]{Sh}). Therefore, Theorem \ref{t1} is more precise
than \cite[Theorem 6]{Sh} (in the commutative case).

Let $f=\sumis p_i \cd m_i$, where $m_1 \stb m_s$ are distinct exponentials
and $p_1 \stb p_s$ are nonzero generalized polynomials. Then the
{\it degree} of $f$ is defined by $\deg f=\sumis (1+\deg p_i )$ if $m_i \ne 1$
for every $i=1\stb s$, and by $\deg f=-1+\sumis (1+\deg p_i )$ if $m_i =1$ for
one of the $i$'s. If $f=0$, then we put $\deg f =-1$.

It is well-known that every $f\in \hbox{\rm GEP}$, $f\ne 0$
has a unique representation of the form $\sumis p_i \cd m_i$, where $m_1 \stb
m_s$ are distinct exponentials and $p_1 \stb p_s$ are nonzero generalized
polynomials. (For Abelian groups this is proved in \cite[Lemma 4.3, p. 41]{Sz}
and in \cite[Lemma 6]{L}. It is easy to check that the proof of
\cite[Lemma 6]{L} works in Abelian semigroups as well.) Thus
$\deg f$ is well-defined for every $f\in \hbox{\rm GEP}$. It is clear that
$\deg f$ equals the usual degree of $f$ if $f\in \hbox{\rm GP}$.

The following result gives the direction (ii)$\akkor$(i) of Theorem \ref{t1}.
\begin{theorem}\label{t2}
Let $f=\sumis p_i \cd m_i$, where $p_i \in \hbox{\rm GP}$, $p_i \ne 0$
$(i=1\stb s)$, and $m_1 \stb m_s$ are distinct exponentials. Put $n_0 =
\max_{1\le i\le s} \deg p_i$ and $k=\deg f$.
\begin{enumerate}[{\rm (i)}]
\item If $m_i \ne 1$ for every $i=1\stb s$, then $f(x_1 +\ldots +x_n )$ is
decomposable of order $k$ for every $n>n_0$.
\item If $m_i =1$ for one of the indices $i=1\stb s$, then
$f(x_1 +\ldots +x_n )$ is decomposable of order $k+1$ for every $n>n_0$.
\end{enumerate}
\end{theorem}
\proof
It is enough to prove the theorem in the special case when 
$s=1$.

The function of $i$ variables $A\colon G^i \to \cc$ is called {\it
$i$-additive,} if 
it is additive in each of its variables (the other variables being fixed).
A function $f_i \colon G\to \cc$ is a {\it monom} of degree $i$, if there
exists a symmetric and $i$-additive function $A_i$ such that
$f_i (x)=A_i (x\stb x)$ for every $x\in G$. 

Let $p\in \hbox{GP}$ be given, and let $\deg p=k$. Then
$p=\sum_{i=0}^k f_i$, where $f_i$ is a monom of degree $i$ for every
$i=1\stb k$, and $f_0$ is constant. (See \cite[Theorem 3]{Dj}.)
Suppose $f_i (x)=A_i (x\stb x)$ $(x\in G)$, where $A_i$ is symmetric and
$i$-additive for every $i=1\stb k$. If $1\le i\le k$, then
\begin{align*}  
f_i (x_1 +& \ldots +x_{k+1} )= A_i (x_1 +\ldots +x_{k+1} \stb x_1 +\ldots +
x_{k+1} )=\cr
  &=\sum_{j_1 \stb j_{i} =1}^{k+1} A_i (x_{j_1}  \stb x_{j_i}) .
\end{align*}
Since $i<k+1$, each term of the sum depends on the variables belonging
to a proper subset of $\{ x_1 \stb x_{k+1} \}$. Grouping the terms appropriately
we find that
$$f_i (x_1 +\ldots +x_{k+1} )=\sum_{j=1}^{k+1} u_{i,j} \quad (i=1\stb k),$$
where $u_{i,j}$ does not depend on $x_j$ $(j=1\stb k+1)$.
Putting $u_j =\sum_{i=1}^{k} u_{i,j} +(f_0 /(k+1))$, we obtain
\begin{equation}\label{e13}
  p (x_1 +\ldots +x_{k+1} )=f_0 +\sumik \sum_{j=1}^{k+1} u_{i,j} =\sum_{j=1}^{k+1}
  u_{j} \cd 1,
\end{equation}
where $u_{j}$ does not depend on $x_j$ $(j=1\stb k+1)$. Therefore,
if $p\in \hbox{GP}$ and $\deg f=k$, then
$p(x_1 +\ldots +x_{k+1} )$ is decomposable of order $k+1$, and then
so is $p(x_1 +\ldots +x_{n})$ for every $n>k$.

Let $f=p\cd m$, where $p\in \hbox{GP}$, $\deg p=k \ge 0$ and $m$ is an
exponential. If $m=1$, then $\deg f=k$, and we obtain that (ii) holds
in this special case.

Suppose $m\ne 1$. Multiplying \eqref{e13} by
$m(x_1 +\ldots +x_{k+1} )=m(x_1 )\cdots m(x_{k+1} )$, we find that
$f(x_1 +\ldots +x_{k+1} )$ is decomposable of order $k+1$, and then
the same is true for $f(x_1 +\ldots +x_{n})$ for every $n>k$. By
$m\ne 1$, we have $\deg f=k+1$, and thus $f$ satisfies (i) (with $k+1$
in place of $k$). \endproof

We note that the bound $k+1$ in (ii) cannot be replaced by $k$, as the following
example shows. Let $F_\om$ denote the free Abelian group generated by
countable infinitely many generators. We shall represent $F_\om$ as
$$F_\om =\{ (t_1 ,t_2 ,\ldots ):  t_i \in \zz \ (i=1,2,\ldots ),
\ \exists \ i_0 , \ t_i =0 \ (i> i_0 ) \} ,$$
where the sum of the elements 
$(t_1 ,t_2 ,\ldots )$ and $(y_1 ,y_2 ,\ldots )$ 
is $(t_1 +y_1 , t_2 +y_2 ,\ldots )$. We take the discrete topology on $F_\om$.

Let $q(x)=\sum_{i=1}^\infty t_i^2$ for every $x=(t_1 ,t_2 ,\ldots ) \in F_\om$.
Then $q \in \hbox{GP}$ and $\deg q=2$. We prove that 
$q(x_1 +\ldots +x_{n})$ is not decomposable of order $2$ for any $n\ge 2$.
Indeed, $q(x_1 +x_2 )$ is not decomposable (of any order), since
otherwise the translates of $q$ would span a linear space of finite dimension,
which is not the case (see \cite{Sz2}). Let $n\ge 3$, and suppose that
$q(x_1 +\ldots +x_{n})$ is decomposable of order $2$. Then we have
\begin{equation*}
q(x_1 +\ldots +x_{n}) =u_1 \cd v_1 +u_2 \cd v_2 ,
\end{equation*}
where the function $u_i$ only depends on the variables belonging to a set
$\emp \ne E_i \subsetneq \{ x_1 \stb x_n \}$, and $v_i$ only depends on
the variables belonging to $\{ x_1 \stb x_n \} \se E_i$ $(i=1,2)$. 
It is easy to check that there are indices $1\le j<k\le n$ such that
$x_j$ and $x_k$ are separated by both $E_1$ and $E_2$; that is,
$x_j \in E_1$ and $x_k \notin E_1$ or the other way around, and
$x_j \in E_2$ and $x_k \notin E_2$ or the other way around.
Fixing $x_i$ at zero for every $i\ne j,k$ we can see that 
$q(x_j +x_k )$ is decomposable of order $2$, which is impossible.

The following result gives the direction (i)$\akkor$(ii) of Theorem \ref{t1}.
\begin{theorem}\label{t3}
If $f\in C(G)$ is such that $f(x_1 +\ldots +x_n )$ is decomposable for some
$n\ge 2$, then $f\in \hbox{\rm GEP}$. More precisely, if $f(x_1 +\ldots +x_n )$
is decomposable of order $k$, then $f\in \hbox{\rm GEP}$, and $\deg f \le k$.
\end{theorem}

As (i) of Theorem \ref{t2} shows, the upper bound $\deg f \le k$ is sharp in
Theorem \ref{t3}.
We prove Theorem \ref{t3} in Section 3, using the results of the next
section.

\section{Lemmas}
In this section we assume that $G$ is a commutative unital
semigroup written additively.
A function $f\colon G\to \cc$ is a {\it polynomial}, if $f=P(a_1 \stb a_n )$,
where $P\in \cc [x_1 \stb x_n ]$ and $a_1 \stb a_n$ are continuous additive
functions mapping $G$ into $\cc$. (These functions are called normal polynomials
in \cite{Sz}.) The functions of the form
$\sumik p_i \cd m_i$, where $p_i$ is a polynomial and $m_i$ is an exponential
for every $i=1\stb k$ are called {\it exponential polynomials.}
If $f\in C(G)$ then we denote by
$V_f$ the linear span of the translates of $f$.

\begin{proposition} \label{p1}
If $f$ is an exponential polynomial, then $\deg f \le \dim V_f$.
\end{proposition}

\proof If $f$ is a polynomial, then even the stronger
statement $\deg f < \dim V_f$ follows from \cite[Lemma 2.8, p. 30]{Sz}.
We give a simple direct proof. We prove by induction on $n=\deg f$. If $n=-1$,
then $f=0$, $V_f =\{ 0\}$, and
$\deg f =-1<0=\dim V_f$. Let $n=\deg f\ge 0$, and suppose that the statement is
true for smaller degrees. There exists an $h\in G$ such that $\deg \De _h f =
n-1$. Fixing such an $h$ we find that $f\notin V_{\De _h f}$, since the elements
of $V_{\De _h f}$ are of degree $\le n-1$. As $V_{\De _h f} +\lan f\ran \su V_f$,
it follows that $\dim V_f >\dim V_{\De _h f}$. By the induction hypothesis we have
$\dim V_{\De _h f}> \deg \De _h f =n-1$, and thus $\dim V_f >n$. This proves
the statement for polynomials.

To prove the general statement, let $f=\sumin p_i \cd m_i$, where
$m_1 \stb$ $m_n$ are distinct exponentials and $p_1 \stb p_n$ are nonzero
polynomials. For every $i$, we have $V_{p_i \cd m_i} =m_i \cd V_{p_i}
=\{ m_i \cd p\colon p\in V_{p_i}\}$, and thus
$$\dim V_{p_i \cd m_i} =\dim V_{p_i} \ge 1+\deg p_i .$$
It is well-known that $p_i \cd m_i \in V_f$ for every $i$. (See, e.g., 
\cite[Lemma 6]{L}. Although the lemma is stated for groups, the proof works in
every Abelian semigroup.) Then we have 
$V_f =V_{p_1 \cd m_1} +\ldots +V_{p_n \cd m_n}$, where the right hand side is a direct sum, as $V_{p_i \cd m_i} \cap \sum_{j\ne i}V_{p_j \cd m_j} =\{0\}$
for every $i=1\stb n$. Therefore, we obtain
$$\dim V_f =\sumin \dim V_{p_i \cd m_i} \ge \sumin (1+\deg p_i ) \ge \deg f.$$
\endproof

We shall need the notion of modified difference operators introduced by
Almira and Sz\'ekelyhidi in
\cite{AS}. If $f,\phi \colon G\to \cc$ and $h\in G$, then we put $\De_{\phi , h}
f(x)=f(x+h)- \phi (h)\cd f(x)$ for every $x\in G$. The statement of the
following lemma is a consequence of \cite[Theorem 17]{AS}
and of \cite[Theorem 8.12, p. 68]{Sz}. In order to make these note as
self-contained as possible, we give an independent direct proof.

\begin{lemma}\label{l1}
\begin{enumerate}[{\rm (i)}]
\item Let $f=\sumik p_i \cd m_i$, where $p_1 \stb p_k$ are generalized
  polynomials and $m_1 \stb m_k$ are exponentials. If $n_i > \deg p_i$
  for every $i=1\stb k$, then we have
\begin{equation}\label{e1}  
  \De_{m_1 , h_1}^{n_1} \ldots \De_{m_k , h_k}^{n_k} f=0
\end{equation}
for every $h_1 \stb h_k \in G$.
\item Let $m_1 \stb m_k$ be different exponentials, and let $n_1 \stb n_k$ be
nonnegative integers. If $f\colon G\to \cc$ is such that \eqref{e1} holds for
every $h_1 \stb h_k \in G$, then $f=\sumik p_i \cd m_i$, where $p_1 \stb p_k$ are
generalized polynomials, and $\deg p_i < n_i$ for every $i=1\stb k$.
\end{enumerate}
\end{lemma}

\proof (i) If $m$ is an exponential, then we have
\begin{align*}
\De _{m,h} \, p_i \cd m_i & = p_i (x+h)\cd m_i (x)\cd m_i (h)-m (h)\cd  p_i (x)
  \cd m_i (x)=\\
  &= [p_i (x+h)\cd m_i (h)-p_i (x) \cd m (h)]\cd m_i (x)=\\
  &=q(x)\cd m_i (x),
\end{align*}
where $q\in \hbox{GP}$, $\deg q\le \deg p_i$, and
$\deg q < \deg p_i$ if $m=m_i$. From this observation the statement follows by
induction on $\sumik \deg p_i$.

\medskip \noindent
(ii) First we prove the statement for $k=1$. We have to prove that if
$m$ is an exponential and $\De_{m,h}^n f=0$ for every $h\in G$, then
$f=p\cd m$, where $p\in \hbox{GP}$ and $\deg p< n$. If $n=0$, then we interpret
the condition as $f=0$. In this case we have $f=0\cd m$, where $\deg 0=-1$, and
so the statement is true. Let $n>0$. It is easy to see that
for every $f\colon G\to \cc$ we have $\De _{m,h}f(x)=m(h)\cd m(x)\cd \De _h (f/m) (x)$. Therefore, by induction we find 
\begin{equation}\label{e2}  
\De _{m,h}^n f=m(h)^n \cd m\cd \De _h^n (f/m) \qquad (h\in G).
\end{equation}
If $\De_{m,h}^n f=0$ for every $h$, then we obtain $\De _h^n (f/m)=0$ for every
$h$. By Djokovi\'c's theorem this condition implies $f/m\in \hbox{GP}$ and
$\deg p\le n-1$ (see \cite[Corollary 1]{Dj}). Then, by $f=(f/m)\cd m$,
we obtain the statement of (ii) for $k=1$.

Suppose that $k>1$ and the statement is true for $k-1$. We argue by induction
on $n=n_1 +\ldots +n_k$. If $n=0$, then $f=0=\sumik 0\cd m_i$, and the statement
is true. Suppose that $n>0$ and that the statement is true for $n-1$.
By symmetry, we may assume $n_1 >0$.
Let $f$ be a function satisfying the condition of (ii).
Let $h_2 \stb h_k \in G$ be fixed, and put
$$f_1 =  \De_{m_2 , h_2}^{n_2} \ldots \De_{m_k , h_k}^{n_k} f.$$
Then $\De_{m_1 , h_1}^{n_1} f_1 =0$ for every $h_1 \in G$.
As we saw above, this implies $f_1 =p\cd m_1$, where $p\in \hbox{GP}$ and
$\deg p <n_1$. Then $\De_{h} \De_{h_1}^{n_1 -1} p =0$ for every $h,h_1 \in G$, and thus
$\De _{m_1 ,h} \De_{m_1 ,h_1}^{n_1 -1} f_1 =0$ for every $h,h_1 \in G$ by \eqref{e2}.

Fix a $h$ such that $m_1 (h)\ne m_2 (h)$. 
Since the difference operators commute, we find that
\begin{equation*}
\De_{m_1 , h_1}^{n_1 -1}  \De_{m_2 , h_2}^{n_2} \ldots \De_{m_k , h_k}^{n_k} \De _{m_1 ,h}f=0
\end{equation*}
for every $h_1 \stb h_k$. By the induction hypothesis, this implies
\begin{equation}\label{e3}  
\De _{m_1 ,h}f=\sumik p_i \cd m_i ,
\end{equation}
where $p_i \in \hbox{GP}$ and $\deg p_i <n_i$ for every $i$.
By the same argument we find that
\begin{equation}\label{e4}  
\De _{m_2 ,h}f=\sumik q_i \cd m_i ,
\end{equation}
where $q_i \in \hbox{GP}$ and $\deg q_i <n_i$ for every $i$.
We have
$$\De _{m_1 ,h}f (x)-\De _{m_2 ,h}f (x)=(m_2 (h)-m_1 (h))\cd f(x),$$
and thus, subtracting \eqref{e4} from \eqref{e3} we get
$$(m_2 (h)-m_1 (h))\cd f =\sumik (p_i -q_i )\cd m_i .$$
Dividing by $(m_2 (h)-m_1 (h))$ we obtain the statement of (ii).
\endproof

\begin{lemma}\label{l2}
Let $V$ be a linear subspace of \hbox{\rm GEP} such that $V\ne \{ 0\}$, and
$\deg f\le N <\infty$
for every $f\in V$. Then there are exponentials $m_1 \stb m_k$ and positive
integers $n_1 \stb n_k$ such that
\begin{enumerate}[{\rm (i)}]
\item $\sumik n_i \le
\begin{cases}
  N   &\text{if} \ 1\notin \{ m_1 \stb m_k \} ,\cr
  N+1 &\text{if}  \ 1\in \{ m_1 \stb m_k \} ;
\end{cases}$
\item every $f\in V$ is of the form $\sumik q_i \cd m_i$, where
$q_i \in \hbox{\rm GP}$ and $\deg q_i < n_i$ $(i=1\stb k)$;
\item there exists a function $f_0 \in V$ such that $f_0 =\sumik p_i \cd m_i$,
where $p_i \in \hbox{\rm GP}$ and $\deg p_i =n_i -1$ for every $i=1\stb k$, and
\item \eqref{e1} holds for every $f\in V$ and for every $h_1 \stb h_k \in G$.
\end{enumerate}
\end{lemma}

\proof
If $f=\sumis p_i \cd m_i$, where $m_1 \stb m_s$ are distinct exponentials
and $p_1 \stb p_s$ are nonzero generalized polynomials, then we define the
{\it spectrum} of $f$ as $\spect f=\{ m_1 \stb m_s \}$. If $f=0$, then we
let $\spect f =\emp$.

Put $S=\bigcup_{f\in V} \spect f$. We show that $S$ has at most $N+1$ elements.
Suppose this is not true. Then there are functions $f_1 \stb f_n \in V$ such
that
\begin{equation}\label{e7}  
\bigcup_{i=1}^n \spect f_i =\{ m_1 \stb m_s \} ,
\end{equation}
where $m_1 \stb m_s$ are distinct exponentials, and $s\ge N+2$. 
Let $f_i =\sum_{j=1}^s p_{i,j} \cd m_j$ $(i=1\stb n)$, where $p_{i,j} \in \hbox{GP}$
for every $i=1\stb n$ and $j=1\stb s$. It follows from \eqref{e7} that for every
$j=1\stb s$ there is an $i$ such that $p_{i,j}\ne 0$. If $c_1 \stb c_n \in \cc$,
then
$$\sumin c_i \cd f_i = \sum_{j=1}^s \left( \sumin c_i \cd p_{i,j} \right) \cd
m_j$$
belongs to $V$. For every $j$, the set of $n$-tuples $(c_1 \stb c_n )$
with $\sumin c_i \cd p_{i,j} =0$ constitutes a linear subspace $L_j$ of
$\cc ^n$. The linear subspace $L_j$ is proper, since $p_{i,j}\ne 0$
for a suitable $i$, and thus $(\de _{1,i}  \stb \de _{n,i} ) \notin L_j$,
where $\de _{k,i}$ is the Kronecker delta. Since $\bigcup_{j=1}^s L_j \subsetneq
\cc ^n$, there is an $n$-tuple $(c_1 \stb c_n )$ such that
$\sumin c_i \cd p_{i,j} \ne 0$ for every $j=1\stb s$, and then the spectrum of
$f=\sumin c_i \cd f_i$ equals $\{ m_1 \stb m_s \}$. Then $\deg f \ge s-1>N$,
which contradicts the condition on $V$.

We proved that $S$ has at most $N+1$ elements. Let $S=\{ m_1 \stb m_k \}$,
where $m_1 \stb m_k$ are distinct exponentials, and $k\le N+1$. Then every
$f\in V$ is of the form $\sumik q_i \cd m_i$, where $q_i =q_i (f) \in
\hbox{GP}$ for every $i=1\stb k$. Since $\deg f\le N$ for every $f\in V$, it
follows that $\deg q_i (f) \le N$ for every $f\in V$ and for every $i=1\stb k$.
Let $n_i =1+\max \{ \deg q_i (f)\colon f\in V\}$ for every $i=1 \stb k$, and
let $f_i \in V$ be such that $\deg q_i (f_i )=n_i -1$ $(i=1\stb k)$.

The set of $k$-tuples $(c_1 \stb c_k )\in \cc ^k$ such that
$\deg \left( \sumjk c_j \cd q_i (f_j )\right) <n_i -1$ is a linear
subspace $M_i$ of $\cc ^k$ for every $i=1\stb k$. We have $M_i \subsetneq
\cc ^k$, as $(\de _{1,i}  \stb \de _{n,i} ) \notin M_i$. Thus $M_i$ is a proper
subspace of $\cc ^k$ for every $i$. Therefore, we have
$\bigcup_{i=1}^k M_i \subsetneq \cc ^k$, and there is an 
$k$-tuple $(c_1 \stb c_k )$ such that
$\deg \left( \sumjk c_j \cd q_i (f_j )\right) =n_i -1$ for every $i=1\stb k$.
Put $f_0 =\sumjk c_j \cd f_j$; then $f_0 \in V$, and (iii) holds.

We have $f_0 =\sumik p_i \cd m_i$, where $p_1 \stb p_k$ are nonzero generalized
polynomials with $\deg p_i =n_i -1$ $(i=1\stb k)$. If $1\notin \{ m_1 \stb
m_k \}$, then $N\ge \deg f_0 =\sumik n_i$. On the other hand, if $1\in \{ m_1
\stb m_k \}$, then $N\ge \deg f_0 =-1+\sumik n_i$, proving (i).
Statement (ii) is obvious from the construction. Finally, (iv) follows
from (ii) and from Lemma \ref{l1}. \endproof

\begin{lemma}\label{l3}
Let $n,N$ be positive integers, and suppose that
\begin{enumerate}[{\rm (i)}]
\item $V$ is a linear subspace of \hbox{\rm GEP}
such that $\deg f\le N$ for every $f\in V$, and
\item $W$ is a linear subspace of $C(G)$ such
that $\dim W \le n$.
\end{enumerate}
Then the number of exponentials contained in $W+V$ is at most $n+N+1$.
\end{lemma}

\proof
By the previous lemma, there are exponentials $m_1 \stb m_k$
such that $k\le N+1$, and every $f\in V$ is of the form $\sumik q_i \cd m_i$,
where $q_1 \stb
q_k \in \hbox{GP}$. We prove that there are at most $n$ exponentials
in $W+V$ different from $m_1 \stb m_k$.

Suppose this is not true, and let $e_1 \stb e_{n+1} \in W+V$ be distinct
exponentials different from $m_1 \stb m_k$. Let $e_i =w_i +v_i$, where
$w_i \in W$ and $v_i \in V$ for every $i=1\stb n+1$. 
Since $\dim W\le n$, the elements $w_1 \stb w_{n+1}$ are linearly
dependent. That is, there are complex numbers $c_1 \stb c_{n+1}$, not all zero,
such that $\sum_{j=1}^{n+1} c_j w_j =0$. Then
$$\sum_{j=1}^{n+1} c_j e_j = \sum_{j=1}^{n+1} c_j v_j \in V.$$ 
Therefore, we have $\sum_{j=1}^{n+1} c_j e_j =\sumik q_i \cd m_i$ with suitable
$q_1 \stb q_k \in \hbox{GP}$. However, the representation of this form is unique, and thus, if $c_j \ne 0$ for a $j$, then $e_j =m_i$ for some $i=1\stb k$,
which is impossible. \endproof

\begin{lemma}\label{l4}
Let $n,N$ be positive integers, and suppose that
\begin{enumerate}[{\rm (i)}]
\item $V$ is a linear subspace of \hbox{\rm GEP}
such that $\deg f\le N$ for every $f\in V$,
\item $W$ is a linear subspace of $C(G)$ such that $\dim W \le n$, and
\item $F$ is a translation invariant linear subspace of $W+V$.
\end{enumerate}
Then $F\su \hbox{\rm GEP}$, and $\deg f\le n+N$ for every
$f\in F$.
\end{lemma}

\proof
Let $m_1 \stb m_k$ and $n_1 \stb n_k$ be as in Lemma \ref{l2}. 
Then we have $\sumik n_i \le N+1$ or $\sumik n_i \le N$
according to whether or not $m_i =1$ for one of the indices $i=1\stb k$.
By Lemma \ref{l3}, the number of exponentials contained in $W+V$ is finite.
Let $m_{k+1} \stb m_{k+s}$ be the
exponentials contained in $W+V$ and different from $m_1 \stb  m_k$.

We show that every function in $F\cap \hbox{GEP}$ is of the form
$\sum_{i=1}^{k+s} q_i \cd m_i$, where $q_i \in \hbox{GP}$ for every $i$.
Suppose $f\in F\cap \hbox{GEP}$ and $f=\sum_{j=1}^t p_j \cd e_j$, where
$e_1 \stb e_t$ are distinct exponentials and $p_1 \stb p_t \in \hbox{GP}
\se \{ 0\}$.
Since $F$ is a translation invariant linear space, it follows that
$e_j \in F\su W+V$ for every $j$. Thus $e_j$ equals one of $m_1 \stb m_{k+s}$
for every $j$, proving the statement.

For every $h_1 \stb h_k \in G$ we put $h=(h_1 \stb h_k )$ and
$$\De _{(h)} = \De_{m_1 , h_1}^{n_1} \ldots \De_{m_k , h_k}^{n_k} .$$
Fix $h_1 \stb h_k$.
By the choice of the numbers $n_i$ and by (iv) of Lemma \ref{l2},
we have $\De _{(h)} v=0$ for every $v\in V$.

It is clear that
$F_{(h)} =\{ \De _{(h)} f \colon f\in F\}$ is a translation invariant linear
subspace of $F$. If $f\in F$, then $f=w+v$ for some $w\in W$ and $v\in V$.
Then $\De _{(h)} f=\De _{(h)} w+\De _{(h)} v= \De _{(h)} w$, and thus
$F_{(h)}$ is a linear subspace of $W_{(h)} =\{ \De _{(h)} w \colon w\in W\}$.
Clearly, $W_{(h)}$ is of dimension $\le n$, and then the same is true for
$F_{(h)}$.

Therefore, by a well-known theorem (see \cite[Theorem 10.1]{Sz}
and the references given on \cite[p. 79]{Sz}), each element of
$F_{(h)}$ is an exponential polynomial. By Proposition \ref{p1},
$\deg \phi \le n$ for every $\phi \in F_{(h)}$.
Therefore, by (i) of Lemma \ref{l1}, we have
$$\De_{m_1 , h_1}^{n+1} \ldots \De_{m_{k+s} , h_{k+s}}^{n+1} \De _{(h)} \phi =0$$
for every $\phi \in F_{(h)}$ and for every $h_{k+1} \stb h_{k+s}\in G$.
According to the definition of $\De _{(h)}$, this implies
$$\De_{m_1 , h_1}^{n_1 +n+1} \ldots \De_{m_k , h_k}^{n_k +n+1}
\De_{m_{k+1} , h_{k+1}}^{n+1} \ldots \De_{m_{k+s} , h_{k+s}}^{n+1} f=0$$
for every $f\in F$ and for every $h_{1} \stb h_{k+s}\in G$. Then, by (ii) of
Lemma \ref{l1}, it follows that $F\su \hbox{\rm GEP}$ and $\sup_{f\in F} \deg f
<\infty$.
We show that $\deg f \le n+N$ for every $f\in F$.

Applying Lemma \ref{l2} with $F$ in place of $V$, we find
a set of indices $J\su \{ 1\stb k+s\}$ and 
positive numbers $d_i$ $(i\in J)$ such that 
every $f\in F$ is of the form $\sum_{i\in J} q_i \cd m_i$, where $q_i
\in \hbox{GP}$ and $\deg q_i < d_i$ $(i\in J)$,
and there is a function $f_0 \in F$ such that $f_0 =\sum_{i\in J} p_i \cd
m_i$, where $p_i \in \hbox{GP}$ and $\deg p_i = d_i -1$
$(i\in J)$. Therefore, we have $\sum_{i\in J} d_i = \deg f_0 +1$ or
$\sum_{i=1}^{k+s} d_i = \deg f_0$
according to whether or not $m_i =1$ for one of the indices $i\in J$.

We put $n_i =0$ for every $i\in J\se \{ 1\stb k\}$. We prove that 
\begin{equation}\label{e12}
\sum_{i\in J} (d_i -n_i ) \le n.
\end{equation}
Suppose this is not true, and let $I$ denote the set of indices
$i\in J$ such that $d_i >n_i$. Then
\begin{equation}\label{e11}
\sum_{i\in I} (d_i -n_i ) \ge \sum_{i\in J} (d_i -n_i ) >n.
\end{equation}
Since $f_0 =\sum_{i\in J} p_i \cd m_i \in F$ and $f$ is a translation invariant
linear space, we have $p_i \cd m_i \in F$ for every $i\in J$. Moreover,
$V_{p_i \cd m_i} =m_i \cd V_{p_i} \su F$ for every $i\in J$. For every
$p\in \hbox{GP}$ we have $\deg \De _h p =(\deg p )-1$ for a suitable $h\in G$,
and thus $V_{p}$ contains generalized polynomials of any degree between $0$
and $\deg p$. Let $p_{i,j} \in V_{p_i}$ be such that 
$\deg p_{i,j} =j$ for every $i\in I$ and $n_i \le j\le d_i -1$.

Let $f_1 \stb f_t$ be an enumeration of the functions
$p_{i,j} \cd m_i$ with $i\in I$ and $n_i \le j\le d_i -1$. Then 
the functions $f_1 \stb f_t$ belong to $F$, and $t=\sum_{i\in I} (d_i -n_i ) >n$
by \eqref{e11}. Since $f_1 \stb f_t \in F$, we have
$f_\nu  =w_\nu  +v_\nu $, where $w_\nu  \in W$, $v_\nu  \in V$ $(\nu =1\stb t)$.
Now $\dim W\le n<t$ implies that there are complex numbers $c_1 \stb c_{t}$,
not all zero, such that $\sum_{\nu =1}^{t} c_\nu  \cd w_\nu  =0$. Let
$f=\sum_{\nu =1}^{t} c_\nu  \cd f_\nu$. Then $f \in F$, and
$f=\sum_{\nu =1}^{t} c_\nu \cd v_\nu\in V$. On the other hand,
$$
f=\sum_{\nu =1}^{t} c_\nu \cd f_\nu   =
\sum_{i\in I} \sum_{j=n_i}^{d_i -1} c_{i,j} \cd p_{i,j} \cd m_i ,$$
where $c_{i,j}$ $(i\in I$, $n_i \le j \le d_i -1)$ is an enumeration of
the numbers $c_\nu$ $(\nu =1 \stb t)$. Thus $f=\sum_{i\in I} q_{i} \cd m_i$,
where $q_i =\sum_{j=n_i}^{d_i -1} c_{i,j} \cd p_{i,j}$ $(i\in I)$. If $i\in I$ 
is such that $c_{i,j}$ is nonzero for at least one $n_i \le j \le d_i -1$,
then $\deg q_i \ge n_i$. If $1\le i\le k$, then $\deg q_i <n_i$ by the choice of
$n_i$, which is impossible. If $i\in J\se \{ 1\stb k\}$,
then $n_i =0$, and thus $\deg q_i \ge n_i$ implies $q_i \ne 0$, contradicting
the choice of $m_1 \stb m_k$. Both cases are impossible, so we must
have \eqref{e12}.

Now we consider three cases. First we assume that $m_i \ne 1$
for every $i\in J\cup \{ 1\stb k\}$. Then we have $\sum_{i=1}^{k} n_i \le N$ and
thus, by \eqref{e12}, we obtain
$$\deg f_0  =\sum_{i\in J} d_i \le n+\sum_{i\in J} n_i \le n+\sum_{i=1}^{k} n_i \le n+N.$$
Next consider the case when $m_i =1$ for some $i \in J$.
Then we have $\sum_{i=1}^{k} n_i \le N+1$ and
thus, by \eqref{e12}, we obtain
$$\deg f_0  =-1 +\sum_{i\in J} d_i \le n-1+\sum_{i\in J} n_i \le n-1+
\sum_{i=1}^{k} n_i \le n+N.$$
Finally, if $m_i =1$ for some $i \in \{ 1\stb k\} \se J$,
Then we have
$$\sum_{i\in J} n_i <\sum_{i=1}^{k} n_i \le N+1$$
and thus, by \eqref{e12}, we obtain
$$\deg f_0  =\sum_{i\in J} d_i \le n +\sum_{i\in J} n_i \le n+N.$$
Summing up: we have $\deg f_0 \le n+N$ in each of the cases. 
Therefore, we have $\deg f\le \deg f_0 \le n+N$ for every $f\in F$. \endproof

\section{Proof of Theorem \ref{t3}}
We follow the argument of the proof of \cite[Theorem 6]{Sh}.
We prove by induction on $n$. If $n=2$ and $f(x_1 +x_2 )$ is decomposable
of order $k$, then $f(x_1 +x_2 )=\sumik u_i (x_1 )\cd v_i (x_2 )$ for every $x_1 ,
x_2 \in G$. This implies that the translates of $f$ belong to the
linear space generated by the functions $u_1 \stb u_k$, and thus $V_f$
is of dimension $\le k$. Therefore, $f$ is an exponential polynomial,
and $\deg f\le k$ by Proposition \ref{p1}.

Let $n>2$, and suppose that the statement is true for $n-1$. Suppose that
$f(x_1 +\ldots +x_n )=\sumik u_i \cd v_i$, where $u_i$ only depends on the
variables belonging to a nonempty set $E_i \subsetneq \{ x_1 \stb x_n \}$,
and $v_i$ only depends on the variables belonging to $\{ x_1 \stb
x_n \} \se E_i$ $(i=1\stb k)$. Switching $u_i$ and $v_i$ if necessary, we may
assume that $x_n \notin E_i$ for every $i$. We may also assume that $E_i =
\{ x_1 \stb x_{n-1}\}$ for every $i\le p$, and $E_i \subsetneq \{ x_1
\stb x_{n-1}\}$ for $p<i\le k$.

If $p=0$, then substituting $x_n =0$ in $f(x_1 +\ldots +x_n )=\sumik u_i \cd v_i$
we can see that $f(x_1 +\ldots +x_{n-1} )$ is decomposable of order $k$. Then,
by the induction hypothesis, we have $f\in \hbox{\rm GEP}$, and $\deg f \le k$,
and we are done. Therefore, we may assume that $p>0$. Note that
the functions $v_1 \stb v_p$ only depend on $x_n$. We have
\begin{equation}\label{e5}  
\sum_{i=1}^p u_i \cd v_i = f(x_1 +\ldots +x_n )-\sum_{j=p+1}^k u_j \cd v_j .
\end{equation}
We may assume that the functions $v_1 \stb v_p$ are linearly independent.
Then there are
elements $g_1 \stb g_p \in G$ such that the determinant $\det |v_i (g_j )|$
is nonzero (see \cite[Lemma 1, p. 229]{AD}). Substituting $x_n =g_\nu$ into
\eqref{e5} for every $\nu=1\stb p$
we obtain a system of linear equations with unknowns $u_1 \stb u_p$.
Since the determinant of the system is nonzero, we find that each of
$u_1 \stb u_p$ is a linear combination of the right hand sides.
None of $u_1 \stb u_k$ depends on $x_n$, so we obtain
\begin{equation}\label{e6}  
u_i =f_i (x_1 +\ldots +x_{n-1} )+\sum_{j=p+1}^k u_j \cd w_{i,j} \qquad (i=1\stb p) ,
\end{equation}
where $f_1 \stb f_p$ are functions of one variable mapping $G$ into $\cc$,
and $w_{i,j}$ only depends on the variables belonging to the nonempty set 
$\{ x_1 \stb x_{n-1}\} \se E_j$ for every $i=1\stb p$ and $j=p+1\stb k$.

Substituting the right hand sides of \eqref{e6} into 
$f(x_1 +\ldots +x_n )=\sumik u_i \cd v_i$ and rearranging the terms we obtain
\begin{equation}\label{e8}
\begin{split}
f(x_1 +& \ldots +x_n ) - \sum_{i=1}^p f_i (x_1 +\ldots +x_{n-1} ) \cd v_i (x_n )=\\
&= \sum_{i=1}^p \sum_{j=p+1}^k  u_j \cd w_{i,j} \cd v_i (x_n ) + \sum_{i=p+1}^k u_j
\cd v_j =\sum_{j=p+1}^k u_j \cd z_j,
\end{split}
\end{equation}
where $z_j =\sum_{i=1}^p (w_{i,j} \cd v_i (x_n ) +v_j)$. Note that,
for every $j=p+1\stb k$, $z_j$ only depends on the 
variables belonging to $\{ x_1 \stb x_n \} \se E_j$, since each function
appearing in the definition of $z_j$ has this property.

Let $\Si$ denote the set of $(k+1)$-tuples
$(\phi , a_1 \stb a_p , b_{p+1}\stb b_k )$ satisfying the following
conditions:
\begin{enumerate}[{\rm (i)}]
\item $\phi \in C(G)$ and $a_1 \stb a_p \colon G\to \cc$,
\item for every $j=p+1\stb k$, $b_j \colon G^n  \to \cc$ is a function only
depending on the variables belonging to $\{ x_1 \stb x_n \} \se E_j$, and 
\item
\begin{equation}\label{e9}
\phi (x_1 +\ldots +x_n ) - \sum_{i=1}^p f_i (x_1 +\ldots +x_{n-1} ) \cd a_i (x_n )=
\sum_{j=p+1}^k u_j \cd b_j
\end{equation} 
for every $x_1 \stb x_n \in G$.
\end{enumerate}
It is clear that $\Si$ is a linear space over $\cc$ under addition and
multiplication by complex numbers coordinate-wise.

Let $F$ denote the set of functions $\phi \in C(G)$ such that
\begin{equation}\label{e10}
(\phi , a_1 \stb a_p , b_{p+1}\stb b_k )\in \Si
\end{equation}
for some $a_1 \stb a_p , b_{p+1}\stb b_k$. By \eqref{e8}, we have $f\in F$.
Since $\Si$ is a linear space over $\cc$, so is $F$.

We show that $F$ is translation invariant. Let $\phi \in F$ and $h\in G$
be arbitrary. Let $(\phi , a_1 \stb a_p , b_{p+1}\stb b_k )\in \Si$.
Then \eqref{e9} holds for every $x_1 \stb x_n \in G$. Replacing $x_n$ by
$x_n +g$ we obtain
$$T_g \phi (x_1 +\ldots +x_n ) - \sum_{i=1}^p f_i (x_1 +\ldots +x_{n-1} ) \cd
T_g a_i (x_n )= \sum_{j=p+1}^k u_j \cd \ol b _j ,$$
where $\ol b _j =T_{(0\stb 0,g)} b_j$. It is clear that $\ol b _j$ only
depends on the variables belonging to $\{ x_1 \stb x_n \} \se E_j$.
Therefore, we have
$$(T_g \phi , T_g a_1 \stb T_g a_p , \ol b _{p+1}\stb \ol b _k )\in \Si ,$$
and $T_g \phi \in F$.

Let $V$ denote the set of functions $\phi -\sum_{i=1}^p a_i (0) \cd f_i$
such that \eqref{e10} holds for suitable
functions $b_{p+1}\stb b_k$. It is clear that $V$ is a linear space.
If $v\in V$, then \eqref{e9} gives
$$v(x_1 +\ldots +x_{n-1} )=\sum_{j=p+1}^k u_j \cd \tilde b _j ,$$
where $\tilde b _j$ is obtained from $b_j$ by putting $x_n =0$. Since
$E_j$ is a proper subset of $\{ x_1 \stb x_{n-1} \}$ for every $j=p+1 \stb k$,
we find that $v(x_1 +\ldots +x_{n-1} )$ is decomposable of order $k-p$.
By the induction hypothesis it follows that $V\su \hbox{GEP}$ and
$\deg v\le k-p$ for every $v\in V$.

Let $W$ denote the linear span of $f_1 \stb f_p$. It is clear that
every $\phi \in F$ is of the form $w+v$, where $w\in W$ and $v\in V$.
Therefore, by Lemma \ref{l4}, we have $\phi \in \hbox{GEP}$ and
$\deg \phi \le p+(k-p)=k$. Since $f\in F$, the proof is complete.
\hfill $\square$

\medskip
\noindent
{\bf Acknowledgement}
The author was supported by the Hungarian National Foundation for Scientific
Research, Grant No. K124749.

\vspace{2cm}

\noindent\textbf{Mikl\'os Laczkovich}\\
E\"otv\"os Lor\'and University\\
Budapest, Hungary\\
{\tt laczk@cs.elte.hu}

\end{document}